\documentclass[letterpaper, 10pt, conference]{ieeeconf}      

\IEEEoverridecommandlockouts                              
\usepackage{graphics} 
\usepackage{epsfig} 
\usepackage{mathptmx} 
\usepackage{times} 
\usepackage{amsmath} 
\usepackage{amssymb}  

\usepackage[mathcal]{euscript} 
\usepackage{mathtools}
\usepackage{mathrsfs}

\usepackage{graphicx}
\graphicspath{{Figures/}} 
\usepackage{caption} 
\usepackage{subfig} 

\usepackage{url}

\usepackage{xcolor}
\definecolor{USred}{rgb}{0.74,0.1,0.1}
\definecolor{USblue}{rgb}{0.2,0.2,0.7}
\definecolor{green1}{cmyk}{0.82,0,1,0.3}

\definecolor{Royalblue}{cmyk}{1,0.30,0.2,0.2}

\newcommand{\Zs}{\mathbb{Z}}
\newcommand{\Es}{\mathbb{E}}
\newcommand{\Rs}{\mathbb{R}}

\newcommand{\boldb}{\mathbf{b}}
\newcommand{\al}[1]{\begin{align} #1\end{align}}
\newcommand{\nn}{\nonumber}
\newcommand*{\qed}{\hfill\ensuremath{\blacksquare}}%

\DeclareMathOperator{\tr}{tr}

\DeclareRobustCommand{\vect}[1]{
	\ifcat#1\relax
	\boldsymbol{#1}
	\else
	\mathbf{#1}
	\fi}

\newtheorem{proposition}{Proposition}

\title{\LARGE \bf A new kernel-based approach for spectral estimation}

\author{Mattia Zorzi 
\thanks{}
\thanks{M. Zorzi is with the Department of Information Engineering, University of Padova, Padova, Italy; email:	 
	   {\tt\small zorzimat@dei.unipd.it}}%
\thanks{}%
}

\begin{document}

\maketitle
\thispagestyle{empty}
\pagestyle{empty}

\begin{abstract}
The paper addresses the problem to estimate the power spectral density of an ARMA zero mean Gaussian process. We propose a kernel based maximum entropy spectral estimator. The latter searches the optimal spectrum over a class of high order autoregressive models while the penalty term induced by the kernel matrix promotes regularity and exponential decay to zero of the impulse response of the corresponding one-step ahead predictor. 
Moreover, the proposed method also provides a minimum phase spectral factor of the process. Numerical experiments showed the effectiveness of the proposed method.
\end{abstract}


\section{Introduction} 
We consider the problem to estimate the power spectral density of an autoregressive and moving average (ARMA) zero mean Gaussian process. The latter has been successfully addressed by the so called ``Tunable High REsolution'' (THREE)-like methods for which a large body of literature is available \cite{ByrnGeoLind2000,GeoLind2003,ferrante2008hellinger,ferrante2012time,zorzi2012estimation,zhu2018existence,Z-14rat,P-F-SIAM-REV} as well as for the multidimensional case \cite{georgiou2006,ringh2018multidimensional,ringh2016multidimensional} and the case of dynamic networks \cite{avventi2012arma,AlpZorzFer2018,ZORZI2019108516,zorzi2019graphical}. According to this approach, the output covariance of a bank of filters is used to extract information on the input power spectral density. More precisely, the estimated spectrum is the one maximizing a suitable entropy-like functional and matching the output covariance matrix. For a comprehensive analysis of the possible entropy-like functionals we refer to \cite{zorzi2014,zorzi15}. In regard to the bank of filters, it has to be designed by the user and it fixes the model class in which we search the optimal spectrum \cite{zorzi2015interpretation}. In general the design of such a bank is very challenging especially in the case that no a priori information about the process is available.

In the simple case that the bank is constituted by $n$ delays, then we obtain the maximum entropy (ME) estimator which is also known as Burg estimator \cite{ulrych1975maximum}.  The latter searches the optimal spectrum over the class of autoregressive (AR) processes of order $n$. The selection of $n$ (which is the unique parameter of the bank of filters) can be performed by resorting to complexity measures such as AIC and BIC, \cite{1100705,schwarz1978estimating}.

Since an (stable) ARMA processes can be approximated by an high order process, one could use the ME spectral estimator 
with $n$ large for ARMA process. However, the resulting estimator will be affected by high variance. In this paper we propose a kernel based ME spectral estimator which searches the optimal spectrum over the class of high order AR processes while the complexity (i.e. degrees of freedom) is controlled through a penalty term induced by kernel matrix. More precisely adopting the Bayesian perspective, the impulse response of the predictor (which is a ``long'' FIR) is modeled as a zero mean Gaussian vector with covariance matrix (i.e. kernel matrix) embedding impulse response regularity and the fact that it decays to zero in an exponential way. Then, we show that the estimated spectrum is given by solving a generalized version of the Yule-Walker equations whose solution leads to a minimum phase spectral factor of the process and thus the estimated spectrum is well defined. Such result is a natural extension of the proof by Stoica and Nehorai \cite{1165162} for the usual Yule-Walker equations. Finally, the kernel matrix is not known: we propose an empirical Bayes approach to estimate it.

It is worth noting that spectral estimation can be performed using the kernel based method proposed in \cite{bottegal2013regularized}. More precisely, it estimates the correlation function, however there is no theoretical guarantee that the corresponding spectrum is positive over the unit circle. An alternative is to consider kernel based prediction error methods (PEM), \cite{NnAaaa,chen2012estimation}. However, it guarantees the predictor stability but it is not guaranteed the system stability \cite{pillonetto2018identification}, i.e. it is not guaranteed that the estimated spectrum is well defined.

The outline of the paper is as follows. In Section \ref{sec:pb_form}
we review the ME spectral estimator. In Section \ref{sec:theta_EST} we  introduce the kernel based ME spectral estimator and in Section \ref{sec:hyper} we propose an empirical Bayes approach to estimate the kernel matrix from the data. In Section \ref{sec:deg} we provide a Fisherian interpretation of the proposed estimator.
In Section \ref{sec:ex} we present some numerical experiments. Finally, in Section \ref{sec:concl} we draw the conclusions.

\section{Maximum Entropy Spectral Estimation}\label{sec:pb_form}
Let $y=\{ y_t,\; t\in\Zs \}$ be a zero mean stationary Gaussian stochastic process. The latter is completely described by the (real-valued) power spectral density \al{ \Phi(e^{j\vartheta})=\sum_{k\in\Zs} r_ke^{j\vartheta k},\;\; \vartheta\in[0,2\pi]}
where $r_k=\Es[y_ty_{t+k}]$ is the $k$-th covariance lag. Notice that $\Phi(e^{j\vartheta})>0$ for any $\vartheta$. We consider the problem to estimate $\Phi$ given a finite length sequence $y^N:=\{ y_1\ldots y_N\}$ extracted from a realization of $y$. Modern high resolution spectral estimators are given by maximizing a suitable entropy-like functional subject to some moment constraints, see \cite{ByrnGeoLind2000}. Here, the model class of the process is characterized by a bank of filters, typically denoted by  $G(z)$, which has to be designed by the user.  However, the choice of $G(z)$ is difficult if no a priori information is available. In the special case that $G(z)=[\,1 \;z^{-1} \ldots\, z^{-n}\,]^T$, with $n\leq N$, and the objective function is the entropy rate, we obtain the ME estimator 
\al{ \label{BURG_ME}\hat \Phi=&\underset{\Phi}{\mathrm{argmax}} \int_{-\pi}^\pi \log \det \Phi (e^{j\vartheta }) \mathrm d\vartheta \nn\\
	&\hbox{ s.t. } \int_{-\pi}^\pi \Phi (e^{j\vartheta }) e^{j\vartheta k}\frac{\mathrm d\vartheta}{2\pi} =\hat r_k, \; \; k=0\ldots n}
where $\hat r_k=\frac{1}{N}\sum_{t=1}^{N-k} y_t y_{t+k}$ is the estimate of $r_k$ using data $y^N$. The latter searches the optimal spectrum over the class of AR models of order $n$ \al{\label{mod_AR}b(z)y_t=e_t}  
where $b(z)=\sum_{k=0}^n b_k z^{-k}$ and $e=\{ e_t,\; t\in\Zs\}$ is white Gaussian noise with unit variance. Indeed $\hat \Phi =(\hat b_{ML}(z)\overline{\hat b_{ML}(z)})^{-1}$ where $\hat b_{ML}(z)$ is characterized as follows. If we define 
\al{\label{def_theta}\boldb=\left[\begin{array}{ccc}b_0 & \ldots & b_n\end{array}\right]^T,}  then $\hat b_{ML}(z)$ corresponds to the maximum likelihood (ML) estimator 
\al{\label{ME_DUAL}\hat \boldb_{ML}=\underset{\boldb\in\Rs^{n+1}}{\mathrm{argmin}} \;\ell(y^N; \boldb) }
where $\ell(y^N)$ is the so called Whittle log-likelihood
\al{\ell(y^N;\boldb):=\frac{N-n}{2}\left(-\log (\boldb^T vv^T \boldb) + \boldb^T \hat \Sigma \boldb\right),}
\al{\label{def_SH}\hat \Sigma =\left[\begin{array}{ccccc} \hat r_0 & \hat r_1 & \ldots & \ldots & \hat r_n\\ \hat r_1 & \hat r_0 &  \ddots & & \vdots \\
\vdots & \ddots & \ddots & & \vdots  \\ \vdots & \ldots &   &   &\vdots  \\
\hat r_n & \ldots & \ldots & \hat r_1 &\hat r_0  \end{array}\right],}
$v=[\,1\; 0\ldots \, 0\,]^T$. Notice that, $\hat \Sigma$ is Toeplitz and positive definite with high probability. It is not difficult to prove that the solution of (\ref{ME_DUAL}) is $\hat \boldb_{ML}=a/\sqrt{a_0}$ where $a=[\, a_0 \; a_1 \ldots \,a_n\,]^T$ is given by solving the Yule-Walker equations  
\al{\label{YW}a=\hat\Sigma^{-1}v.} It is well known that $\hat b_{ML}(z)$ has zeros inside the unit circle and thus $\hat b_{ML}(z)^{-1}$ is a minimum phase spectral factor of $\hat \Phi$.  Moreover, $\hat \Phi>0$, i.e. the spectrum is well defined. Although in (\ref{BURG_ME}) the bank of filters is very simple, we need to choose $n$, that is the order of the AR model. The latter step is typically addressed by resorting to complexity measures such as AIC and  BIC \cite{1100705,schwarz1978estimating}. However, the selection of $n$ is not trivial and may compromise the optimality properties, in particular for shorter data records; moreover, it corresponds to a combinatorial optimization problem. In what follows we propose a kernel based maximum entropy estimator which searches the optimal over a large model class (i.e. high order AR models) while the complexity is controlled by a regularization term induced by the kernel matrix.

\section{Spectral estimation using a Gaussian prior}\label{sec:theta_EST}

We address the problem to estimate the power spectral density $\Phi$ of a zero mean ARMA Gaussian stationary process $y$. Let $w(z)$ be a minimum phase rational spectral factor of $\Phi$, then 
\al{w(z)^{-1}y_t=e_t}
where $e_t$ is white Gaussian noise with variance equal to one. Since $w(z)^{-1}$ is rational, we can rewrite it using the Laurent expansion:
\al{w(z)^{-1}=\sum_{k=0}^\infty b_k z^{-k}.} Notice that $-b_0^{-1}b_k$, with $k\geq 1$, is the impulse response of the one-step ahead predictor (see equation (\ref{eq_pred}) below). Since $w(z)$ is minimum phase, we have that the sequence $b_k$, with $k\geq 0$, belongs to $\ell_1$, i.e. $\sum_{k=0}^\infty |b_k|<\infty$. This implies that 
\al{w(z)^{-1}\approx b(z):=\sum_{k=0}^n b_kz^{-k}}
where $n$ is sufficiently large. In other words, we can always approximate an ARMA process with an high order AR process. Let $\boldb$ be defined as in (\ref{def_theta}). Given the dataset $y^N$, then we could estimate it by solving (\ref{ME_DUAL}). However, the corresponding estimator will be affected by high variance, i.e. the complexity of the model is not well controlled.

In order to counteract this aspect, we adopt the Bayesian perspective. More precisely, we model $\boldb$ as a Gaussian random vector with zero mean and covariance matrix $\lambda K$. $\lambda>0$ is called scale factor and $K>0$ is called kernel matrix. Therefore, its probability density is:
\al{p(\boldb)=\frac{1}{\sqrt{(2\pi)^{n+1} \det K}} \exp\left(-\frac{1}{2\lambda} \boldb^TK^{-1}\boldb\right).} 
The, the negative log-likelihood of $y^N$ and $\boldb$ is
\al{\ell&(y^N,\boldb)=-\log p(y^N,\boldb)=-\log p(y^N|\boldb)-\log p(\boldb)\nn\\
& = \frac{N-n}{2}\left(-\log (\boldb^T vv^T \boldb) +  \boldb^T \hat \Sigma \boldb\right)\nn\\ & \hspace{0.4cm}+\frac{1}{2}\log \det (\lambda K) +\frac{1}{2\lambda }\boldb^T K^{-1}\boldb +c_1  }
where $c_1$ is a constant term not depending on $\boldb$ and $K$.
Assume that we have a preliminary estimate of $b_0$, say $\tilde b_0$, then 
\al{\log & (\boldb^T v v^T \boldb )=\log b_0^2=2\log b_0 \approx 2\log \tilde b_0+\frac{2}{\tilde b_0}(b_0-\tilde b_0)\nn\\ &=   2\log \tilde b_0+\frac{2}{\tilde b_0}(\boldb^T v-\tilde b_0)}
accordingly
\al{\label{nlk1}\ell&(y^N,\boldb) \approx \frac{N-n}{2}\left(- \frac{2}{\tilde b_0}\boldb^T v +\boldb^T \hat \Sigma \boldb\right)\nn\\ & \hspace{0.4cm}+\frac{1}{2}\log \det (\lambda K )+\frac{1}{2\lambda}\boldb^T K^{-1}\boldb +c_2   }
where $c_2$ is a constant term not depending on $\boldb$ and $K$. It is worth noting that we can rewrite the AR model in (\ref{mod_AR}) as
\al{\label{eq_pred}y_t=a(z)y_t+u_t}
where $a(z)=1-b_0^{-1}b(z)$ and $u_t$ is white noise with variance equal to $b_0^{-2}$. As suggested in \cite{148344}, a simple an effective way to estimate $b_0^{-2}$ is to consider a low order AR model and solve (\ref{BURG_ME}). Then, notice that adding the constant term $(N-n) v^T \hat\Sigma^{-1}v/2\tilde b_0^2$ in (\ref{nlk1}) we obtain
\al{\ell&(y^N,\boldb) \approx \frac{1}{2}\| \tilde v-\Phi \boldb\|^2+\frac{1}{2}\log \det (\lambda K) +\frac{1}{2\lambda }\boldb^T K^{-1}\boldb +c_3 \nn  }
where
\al{\label{defvPhi}\tilde v= \frac{\sqrt{N-n}}{\tilde b_0} L^{-1} v, \; \; \Phi=\sqrt{N-n} L^T, \; \; \hat\Sigma=LL^T }
 $c_3$ is a constant term not depending on $\boldb$ and $K$.
Notice that $L$ is an arbitrary square root decomposition of $\hat \Sigma$, e.g. the Cholesky decomposition.  
Accordingly, the ML estimator for $\boldb$ is 
\al{\label{ReLS}\hat \boldb_{ML}=\underset{\boldb\in\Rs^{n+1}}{\mathrm{argmin}}\; \ell(y^N,\boldb)=\underset{\boldb\in\Rs^{n+1}}{\mathrm{argmin}} \| \tilde v-\Phi \boldb\|^2 +\lambda^{-1} \|\boldb	\|^2_{ K^{-1}}. }
Therefore, $\hat \boldb_{ML}$ is the solution to a regularized Least squares problem. In particular, it is not difficult to prove that 
\al{\label{th_opt}\hat \boldb_{ML}&=\lambda K\Phi^T (	\lambda \Phi K \Phi^T +I)^{-1}\tilde v\nn\\
& =\hat b_0^{-1}\left(\hat \Sigma+((N-n) \lambda  K)^{-1}\right)^{-1}v.}  Notice that $\hat \boldb_{ML}$ does not depends on the particular decomposition $\hat \Sigma=LL^T$.
It is worth noting that the role of the penalty term $\lambda^{-1}\|\boldb\|^2_{K^{-1}}$ is to reduce the variance of the estimator since $n$ is chosen large. In case that $\lambda\rightarrow \infty$, i.e. there is no a priori information on $\boldb$, then we obtain the usual maximum entropy problem in (\ref{ME_DUAL}).

It remains to design the kernel matrix $K$ that encodes the a priori information on $b(z)$. For instance, we want that $\hat \boldb_{ML}$ corresponds to a polynomial $\hat b_{ML}(z)$ whose zeros are inside the unit circle. Moreover, $b_k$ should decay to zero in an exponential way as $k$ increases because the latter are the coefficients of the Laurent expansion of a transfer function  whose poles are inside the unit circle. In what follows we propose two kernel matrices satisfying such requirements. The latter has been already used in kernel based PEM methods in a successful way, see e.g.  \cite{NnAaaa,chen2012estimation}.
 
\subsection{Diagonal (DI) kernel}
\begin{figure}
\includegraphics[height=0.27\textwidth,width=0.5\textwidth]{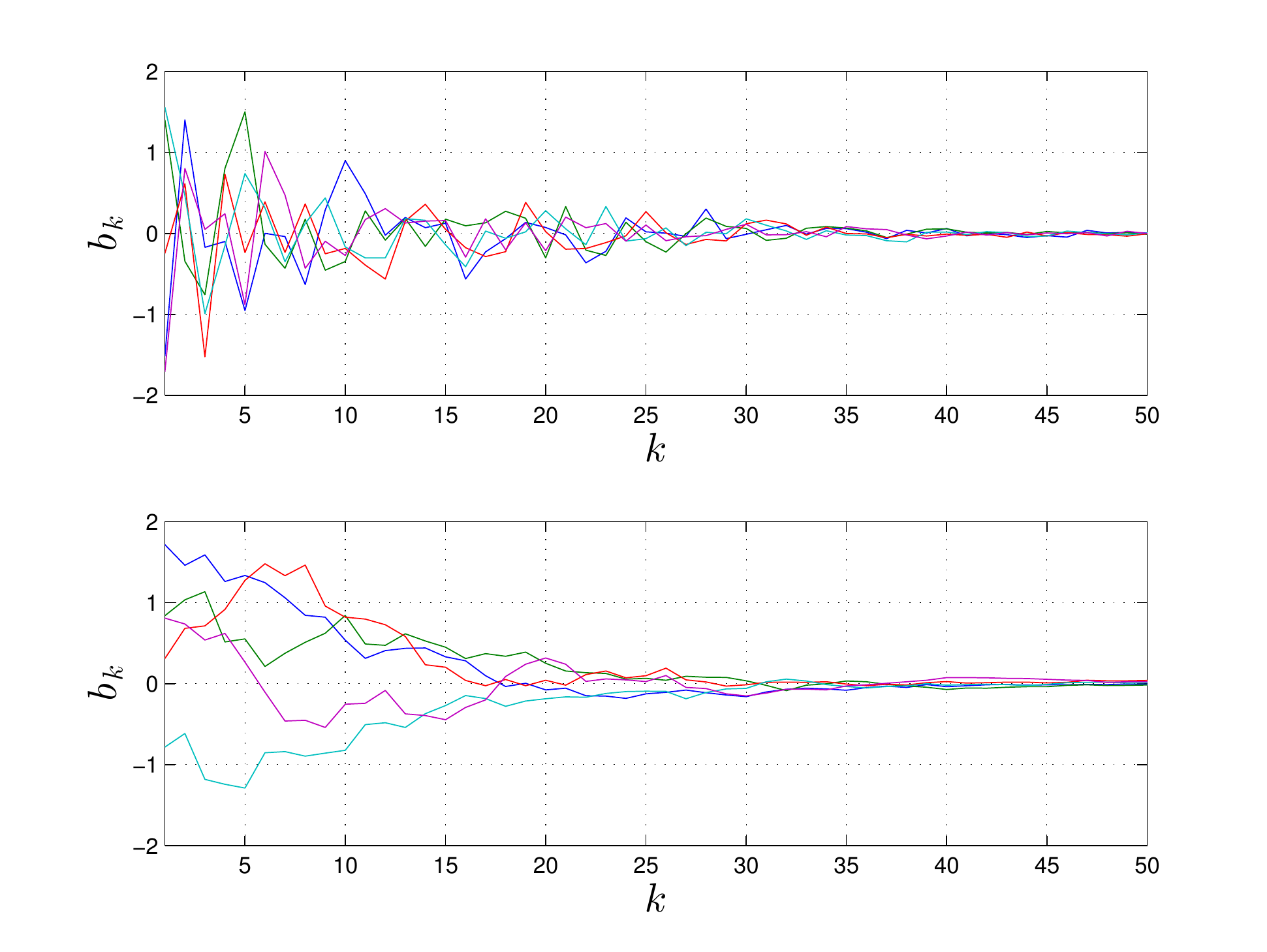}
\caption{{\em Top panel.} Five realizations of $\boldb$ using the Diagonal kernel. {\em Bottom panel.} Five realizations of $\boldb$ using the Tuned-Correlated  kernel. }\label{fig_rel}
\end{figure}

We consider the kernel matrix defined as 
\al{K_{DI}=\mathrm{diag}(\beta,\beta^2,\beta^3\ldots \beta^{n+1})  }
where $\beta\in(0,1)$. In Figure \ref{fig_rel} (top panel) we show five realizations of $\boldb \sim\mathcal N (0,\lambda K_{DI})$ with $\lambda=1$ and $\beta=0.85$. 
 As we see, in all these realizations $b_k$ decays in an exponential way. Indeed, for any realization of $\boldb$ we have that $\|\boldb\|^{2}_{K^{-1}}<\infty$ almost surely which implies that $b_k$ decays in an exponential way almost surely.  

\begin{proposition} \label{prop1}Let $\hat \boldb_{ML}$ be the estimate using $K_{DI}$. If the corresponding polynomial $\hat b_{ML}(z)$ is nonnull, then it has zeros inside the unit circle. 
	\end{proposition} \proof By (\ref{th_opt}) we have that 
	\al{\label{rel2}(\hat \Sigma+R) \hat \boldb_{ML}=\tilde b_0^{-1} v} where \al{R=((N-n)\lambda K)^{-1}:=\mathrm{diag} (r_1,r_2\ldots r_{n+1})}
	and $r_{k}=((N-n)\lambda\beta^{k})^{-1}$. Thus, $r_{k+1}\geq r_k>0$. Let $\breve z$ be a zero of $\hat b_{ML}(z)$. If $\breve z=0$, then it is inside the unit circle and we have finished. 
Otherwise we have $\breve z\neq 0$. Define  
\al{u_k=\left\{\begin{array}{ll} \hat b_0, & k=1\\
 \breve z u_{k-1}+\hat b_{k-1}, & k=2\ldots n\\
\end{array}\right.} where $\hat \boldb_{ML}=[\,\hat b_0 \; \hat b_1\ldots \, \hat b_n\,]^T$. Accordingly, \al{\label{rel1} u=\breve z \tilde u+\hat \boldb_{ML}}
where $u=[\,u_1 \; u_2\ldots \, u_n\; 0\,]^T$ and $\tilde u=[\,0\; u_1\; u_2 \ldots \, u_n\,]^T$. Using (\ref{rel1}), we have 
\al{u^* & (\hat \Sigma+R) u=(\overline{\breve z} \tilde u^*+\hat\boldb_{ML}^T)(\hat \Sigma+R)( \breve z  \tilde u +\hat\boldb_{ML} )\nn\\ &=|\breve z |^2 \tilde u^*(\hat \Sigma+R) \tilde u +\hat\boldb_{ML}^T(\hat \Sigma+R)\hat\boldb_{ML}+2\Re[\breve z  \hat\boldb_{ML}^T(\hat \Sigma+R) \tilde u ]\nn\\
&=|\breve z |^2 \tilde u^*(\hat \Sigma+R) \tilde u +\hat\boldb_{ML}^T(\hat \Sigma+R)\hat\boldb_{ML}+2\tilde b_0^{-1}\Re[\breve z  v^T\tilde u ]}
where in the last equality we exploited (\ref{rel2}). Notice that $v^T\tilde u =0$, accordingly the last term on the right hand side is equal to zero. Moreover $u^* \hat \Sigma u=\tilde u^* \hat \Sigma \tilde u$ because $\hat \Sigma$ is a Toeplitz matrix. Therefore,
\al{\label{rel3}u^* & (\hat \Sigma+R) u=|\breve z |^2  u^* \hat \Sigma  u +|\breve z |^2 \tilde u^* R \tilde u  +\hat\boldb_{ML}^T(\hat \Sigma+R)\hat\boldb_{ML}.} Finally notice that 
\al{\tilde u^* R\tilde u=\sum_{k=1}^n r_{k+1} |u_k|^2\geq \sum_{k=1}^n r_{k} |u_k|^2=u^* R u.} Using the latter inequality in (\ref{rel3}), we obtain
\al{\label{rel4}u^* & (\hat \Sigma+R) u\geq |\breve z |^2  u^* (\hat \Sigma+R)  u +\hat\boldb_{ML}^T(\hat \Sigma+R)\hat\boldb_{ML}} and thus 
\al{\label{rel5}(1-|\breve z |^2)u^* & (\hat \Sigma+R) u\geq    \hat\boldb_{ML}^T(\hat \Sigma+R)\hat\boldb_{ML}.}  The right hand side in (\ref{rel5}) is positive because $\hat \Sigma+R>0$ and $\hat \boldb_{ML}\neq 0$ (otherwise $\hat b_{ML}(z)=0$). Moreover,  
$u^*  (\hat \Sigma+R) u>0$ because $u\neq 0 $ (otherwise $\hat \boldb_{ML}=0$). Accordingly, we have $1-|\breve z|^2>0$, that is $\breve z$ is inside the unit circle. \qed

It is worth noting the proof above is in the same spirit of the one in \cite{1165162} for the Yule-Walker equations in (\ref{YW}). 

\subsection{Tuned-Correlated (TC) kernel}
We consider the kernel matrix $K_{TC}$ whose entry in position $(t,s)$ is defined as 
\al{\label{def_TC}[K_{TC}]_{ts}=\beta^{\max(t,s)}-\beta^{n+2}, \; \; t,s=1\ldots n+1}
where $\beta\in(0,1)$. It is worth noting that $K_{TC}$ in (\ref{def_TC}) is a modified version of the one proposed in \cite{chen2012estimation}. Indeed, the latter is defined as \al{[\tilde K_{TC}]_{t,s}=\beta^{max(t,s)},\; \; t,s=1\ldots n+1.}
However, for $n$ large the two kernel matrices $\tilde K_{TC}$ and $K_{TC}$ are very similar. 
 In Figure \ref{fig_rel} (bottom panel) we show five realizations of $\boldb\sim  \mathcal N (0,\lambda K_{DI})$ with $\lambda=1$ and $\beta=0.85$. As we see, in all these realizations $b_k$ decays in an exponential way. Moreover, these  impulse responses are more smooth than the ones using the diagonal kernel. Also in this case, for any realization of $\boldb$ we have that $b_k$ decays in an exponential way almost surely.  It is not difficult to prove  that 
\al{K_{TC}=(FDF^T)^{-1}}
 where
\al{\label{FD}F &=\left[\begin{array}{ccccc} 1& 0 & \ldots & \ldots & 0\\ 
-1 & 1 &  \ddots & & \vdots \\
0 & \ddots & \ddots & & \vdots  \\ \vdots & \ldots &   &   &\vdots  \\
0 & \ldots & \ldots & -1 & 1  \end{array}\right] \\
D&= \frac{1}{\beta-\beta^2}\mathrm{diag}\left(1,\frac{1}{\beta},\frac{1}{\beta^{2}}\ldots \, \frac{1}{\beta^{n-1}},\frac{1}{\beta^n}\right).}

\begin{proposition} Let $\hat \boldb_{ML}$ be the estimate using $K_{TC}$. If the corresponding polynomial $\hat b_{ML}(z)$ is nonnull, then it has zeros inside the unit circle. 
	\end{proposition} \proof Equation (\ref{rel2}) still holds  with 
	$R=F\tilde DF^T$ where \al{\tilde D=\mathrm{diag} (d_1,d_2\ldots d_{n+1})}
	and $d_{k}=((N-n)\lambda\beta^{k-1}(\beta-\beta^2))^{-1}$ for $k=1\ldots n+1$. Thus, $d_{k+1}\geq d_k>0$ for $k=1\ldots n$. 
Let $\breve z$ a zero of $\hat b_{ML}(z)$. If $\breve z=0$, then it is inside the unit circle and we have finished. 
Otherwise we have $\breve z\neq 0$. We define $u$ and $\tilde u $ as in the proof of Proposition \ref{prop1}. Accordingly, (\ref{rel1}) holds. Using the same arguments in the proof of Proposition \ref{prop1}, we obtain (\ref{rel3}). Then, we have  
\al{\tilde u^* R\tilde u&= d_1 |u_1|^2 +\sum_{k=1}^{n} d_{k+1} |u_k-u_{k+1}|^2  +d_{n+1}|u_n|^2\nn\\ &\geq d_1 |u_1|^2 +\sum_{k=1}^{n} d_{k} |u_k-u_{k+1}|^2+d_{n+1}|u_n|^2 \nn\\
  &\geq \sum_{k=1}^{n} d_{k } |u_k-u_{k+1}|^2  = u^*  R u\nn }  and thus $\tilde u^*  R\tilde u\geq u^*  R u$. Accordingly, using the same arguments in the proof of Proposition \ref{prop1} we obtain (\ref{rel5}) and thus  $\breve z$ is inside the unit circle. \qed

\section{Hyperparameters estimation}\label{sec:hyper} 

The maximum likelihood estimator of Section \ref{sec:theta_EST} depends on $\eta:=[\, \lambda\; \beta\,]^T$ which are called hyperparameters in the machine learning community. $\eta$ is not known and has to be estimated form the data $y^N$. According to the empirical Bayes perspective, \cite{friedman2001elements}, an estimate of $\eta$ is given by maximizing the probability density of $y^N$ under the model AR of order $n$ in (\ref{mod_AR}) with $\boldb\sim  \mathcal N(0,\lambda K_\beta)$ where we made explicit the fact that $K$ depends on $\beta$.   
In our case, such probability density is 
\al{p(y^N)=c_4\int_{\Rs^{n+1}} \exp(-\ell (y^N,\boldb)) \mathrm d \boldb } 
where $c_4$ is a positive constant term not depending on $\boldb$ and $\eta$. Accordingly, an estimator for $\eta$ is 
\al{\hat \eta=&\underset{\eta}{\mathrm{argmax}}\, \ell(y^N) \nn\\
& \hbox{ s.t. } \lambda>0,\; \; 0<\beta<1}
where $\ell(y^N)=-\log p(y^N)$ is the so called negative log-marginal likelihood function. In order to find an analytical expression for $\ell(y^N)$ we need the following results. 
  
\begin{proposition} \label{prop_int}Let $h\,:\, \Rs^{n+1}\longrightarrow\Rs$ be such that its Hessian matrix $\partial^2 h(\boldb)/ \partial \boldb^2=H$ is constant and positive definite.Then, \al{\int_{\Rs^{n+1}} \exp(-h(\boldb)) \mathrm d \boldb =\frac{1}{\sqrt{\det(2\pi H)}}\exp(-h( \boldb^\circ))} where $\boldb^\circ$ minimizes $h(\boldb)$, see \cite{de1981asymptotic}.
\end{proposition}

In our case $h(\boldb)=\ell(y^N,\boldb)$ and  
\al{H=\frac{\partial^2\, \ell(y^N,\boldb)}{\partial \,\boldb^2}= \Phi^T\Phi+\lambda^{-1} K_\beta^{-1}} 
which is constant and positive definite. Moreover, $\boldb^\circ=\hat \boldb_{ML}$ and it is not difficult to see that
\al{\ell(y^N,\hat \boldb_{ML})= \frac{1}{2}\tilde v^T (\lambda \Phi K_\beta\Phi^T +I)^{-1}\tilde v+\frac{1}{2} \log\det(\lambda K_\beta) +c_2.}
Then, we can apply Proposition \ref{prop_int}, thus
\al{p(y^N)= \frac{1}{\sqrt{\det(2\pi (\Phi^T\Phi +\lambda^{-1}K_\beta^{-1})) }}\exp(-\ell(y^N,\hat \boldb_{ML})) } 
and \al{\label{an_ell}\ell(&y^N) = \frac{1}{2} \log \det(\Phi^T\Phi +\lambda^{-1}K_\beta^{-1})+ \frac{1}{2}\tilde v^T (\lambda \Phi K_\beta\Phi^T +I)^{-1}\tilde v\nn\\ &\hspace{0.5cm}+\frac{1}{2} \log\det(\lambda K_\beta) +c_5\nn\\
 &= \frac{1}{2} \log \det(\lambda \Phi K_\beta\Phi^T+I)+ \frac{1}{2}\tilde v^T (\lambda \Phi K_\beta\Phi^T +I)^{-1}\tilde v+c_5
}
where $c_5$ is a constant term not depending on $\eta$.

The resulting spectral estimation  procedure can be summarized in the following steps:
\begin{enumerate}
\item Compute a preliminary estimate $\tilde b_0$ of $b_0$ from $y^N$ using a low order AR model 
\item Set $n$ large and such that $n\leq N$
\item Compute the Toeplitz matrix $\hat \Sigma$ in (\ref{def_SH}) from $y^N$  
\item Compute a square root of $\hat \Sigma$
\item Compute $\Phi$ and $\tilde v$ as in (\ref{defvPhi})
\item Compute $\hat \eta=[\,\hat \lambda \; \hat \beta\, ]^T$ which minimizes (\ref{an_ell})  
\item $\hat \boldb_{ML}=\hat \lambda K_{\hat \beta} \Phi^T (\hat \lambda \Phi K_{\hat \beta}\Phi^T+I)^{-1}\tilde v$.
\end{enumerate}
It is worth noting that (\ref{an_ell}) is not a convex function. However, there is a large body of literature dealing with the global minimization problem of functions of type (\ref{an_ell}) under constraints $\lambda>0$ and $0<\beta<1$, se e.g. \cite{MULTILPLE_TAC}.

\section{Degrees of Freedom}\label{sec:deg}
Problem (\ref{ReLS}) has also an interpretation in the Fisherian perspective, i.e. when $\boldb$ is understood as an unknown but fixed vector. In this case, (\ref{ReLS}) is a regularized ML estimator of $\boldb$. Since (\ref{ReLS}) is a regularized Least squares problem, the corresponding degrees of freedom of the estimated model are defined as  \cite{tibshirani2015degrees}:
\al{df(\eta)&=\tr[\Phi (\Phi^T\Phi+(\lambda K_\beta)^{-1})^{-1}\Phi^T]\nn\\
&=\tr[ (\hat \Sigma+((N-n)\lambda K_\beta)^{-1})^{-1}\hat \Sigma].}
It is well known that $df(\eta)\in\Rs$ (and nonnegative) measures the complexity of the model: the higher $df(\eta)$ is, the more complex the model is. In the special case that $\lambda =\infty$, i.e. no regularization is present, we have:
	\al{ df(\eta)=\tr[\hat \Sigma^{-1} \hat \Sigma]=n+1}
that is, as expected, the complexity is an integer number corresponding to the number of parameters to estimate for the AR model. We conclude that regularization allows to control the complexity, through $\eta$, in a finer way than the case without regularization. 

\section{Example}\label{sec:ex}
We test the performance of the proposed kernel based spectral estimator using the DI and TC kernel with the ME spectral estimator. An alternative way to estimate the spectrum is to use a kernel-based PEM method introduced in \cite{NnAaaa}. The latter estimates an high order AR model by minimizing the one-step ahead prediction error of the model 
plus the penalty term induced by the DI or the TC kernel. Note that this method does not guarantees that the $b(z)$ has zeros inside the unit circle, and thus that the estimated power spectral density is well defined.
In what follows we will use the following notation for the spectral estimators:
\begin{itemize}
\item ME denotes the ME estimator equipped with the BIC criterium for selecting the optimal $n$ over the set $\{1,2\ldots \, 50\}$
\item ME-DI denotes the kernel based ME estimator using the DI kernel
\item ME-TC denotes the kernel based ME estimator using the TC  kernel
\item PEM-DI denotes the kernel based PEM estimator using the  DI kernel
\item PEM-TC denotes the kernel based PEM estimator using the TC kernel.
\end{itemize} For all the kernel based methods we set $n=50$.

\begin{figure}
\includegraphics[width=0.5\textwidth]{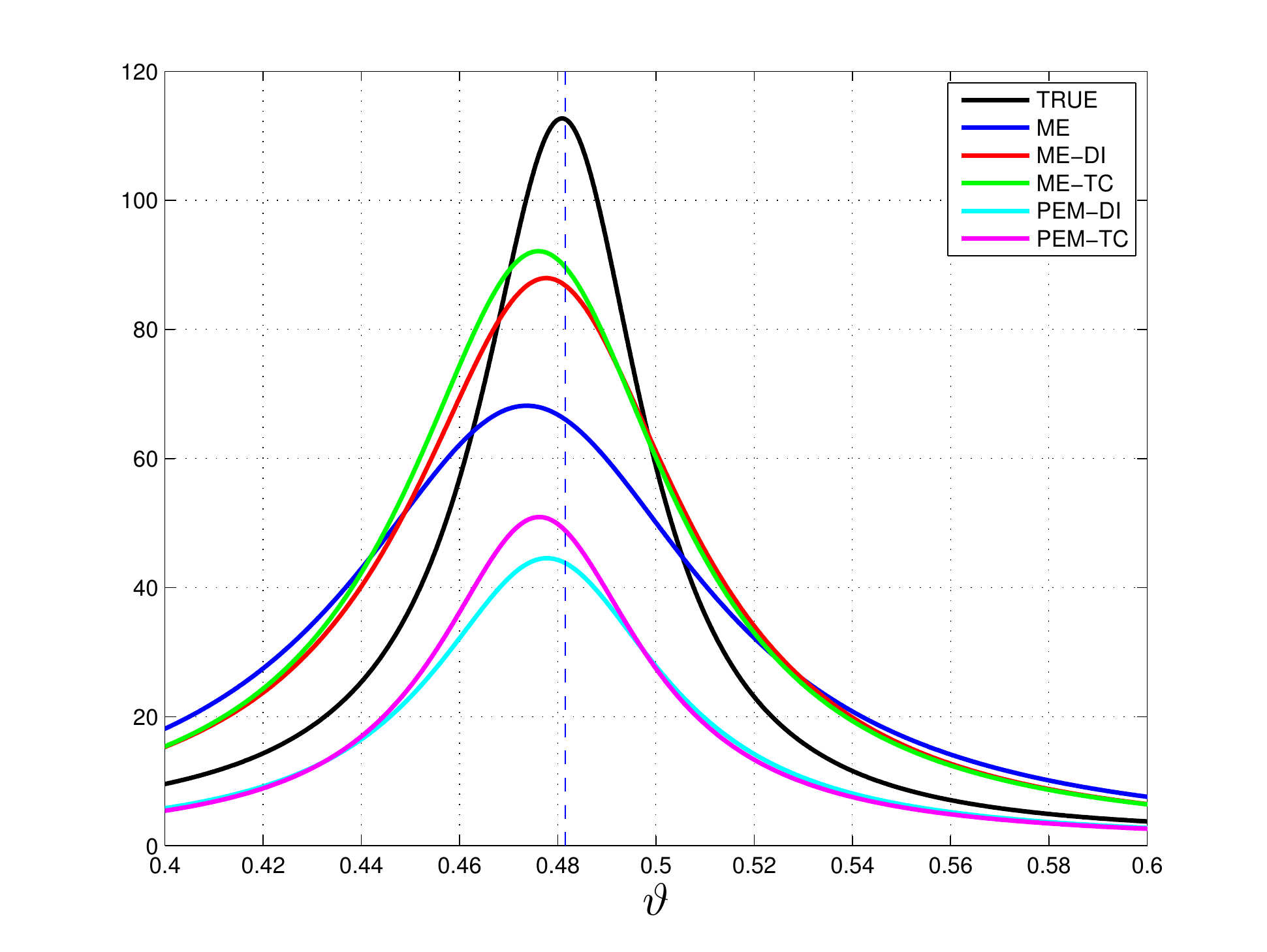}
\caption{Estimated spectra versus the true spectrum over the frequency interval $[0.4,0.6]$. The dashed line depicts the frequency of the poles of  the true model.}\label{fig_trial}
\end{figure}

In the first experiment we generate a dataset $y^N$ with $N=500$ 
from the ARMA process 
\al{\label{1st_exp}y_t=\sqrt{2}\frac{(z-z_0)(z-\overline z_0 ) }{(z-p_0)(z-\overline p_0 )}e_t}
where $z_0=0.85e^{j0.52}$ and  $p_0=0.98e^{j0.482}$. 
Each estimator returns a polynomial $\hat b(z)$, then the corresponding spectral estimator is $\hat \Phi(z)=(\hat b(z) \hat b(z^{-1}) )^{-1}$. 
Figure \ref{fig_trial} depicts the estimated spectra versus the true spectrum.
 As we can see ME-DI and ME-TC performs better than ME. Finally, the PEM method provides the worst performance. Regarding the degrees of freedom we have: $df_{ME}=\hat n_{BIC}+1=8$ $df_{ME-DI}(\hat \eta)=10,8131$ and  $df_{ME-TC}(\hat \eta)=8,8943$ where $\hat n_{BIC}$ is the optimal AR order according to the BIC criterium.

In the second experiment we consider a Monte Carlo study composed by $100$ runs. In each run we generate a dataset $y^N$ with $N=500$ of an ARMA process $y_t=w(z) e_t$ where $w(z)$ is a minimum phase spectral factor with 3 zeros and poles. Each poles has modulus equal to $0.98$ while  the phase is drawn form a uniform distribution over $[0,\pi]$. Each zero has modulus equal to 0.85 while the phase is chosen randomly (according to a uniform distribution) but sufficiently close to the one of a pole. More precisely, if $\vartheta_z$ denotes the phase of the zero, then there exists a pole with phase $\vartheta_p$ such that $|\vartheta_z-\vartheta_p|\leq 0.06$. Then, for each estimator we evaluate the reconstruction error which is defined as:
\al{e=\frac{\int_{-\pi}^\pi \| \hat \Phi(e^{j\vartheta})-   \Phi(e^{j\vartheta}) \|^2\mathrm d \vartheta}{\int_{-\pi}^\pi \|   \Phi(e^{j\vartheta}) \|^2\mathrm d \vartheta}}
where $\Phi$ denotes the true spectrum.
\begin{figure}
\includegraphics[width=0.5\textwidth]{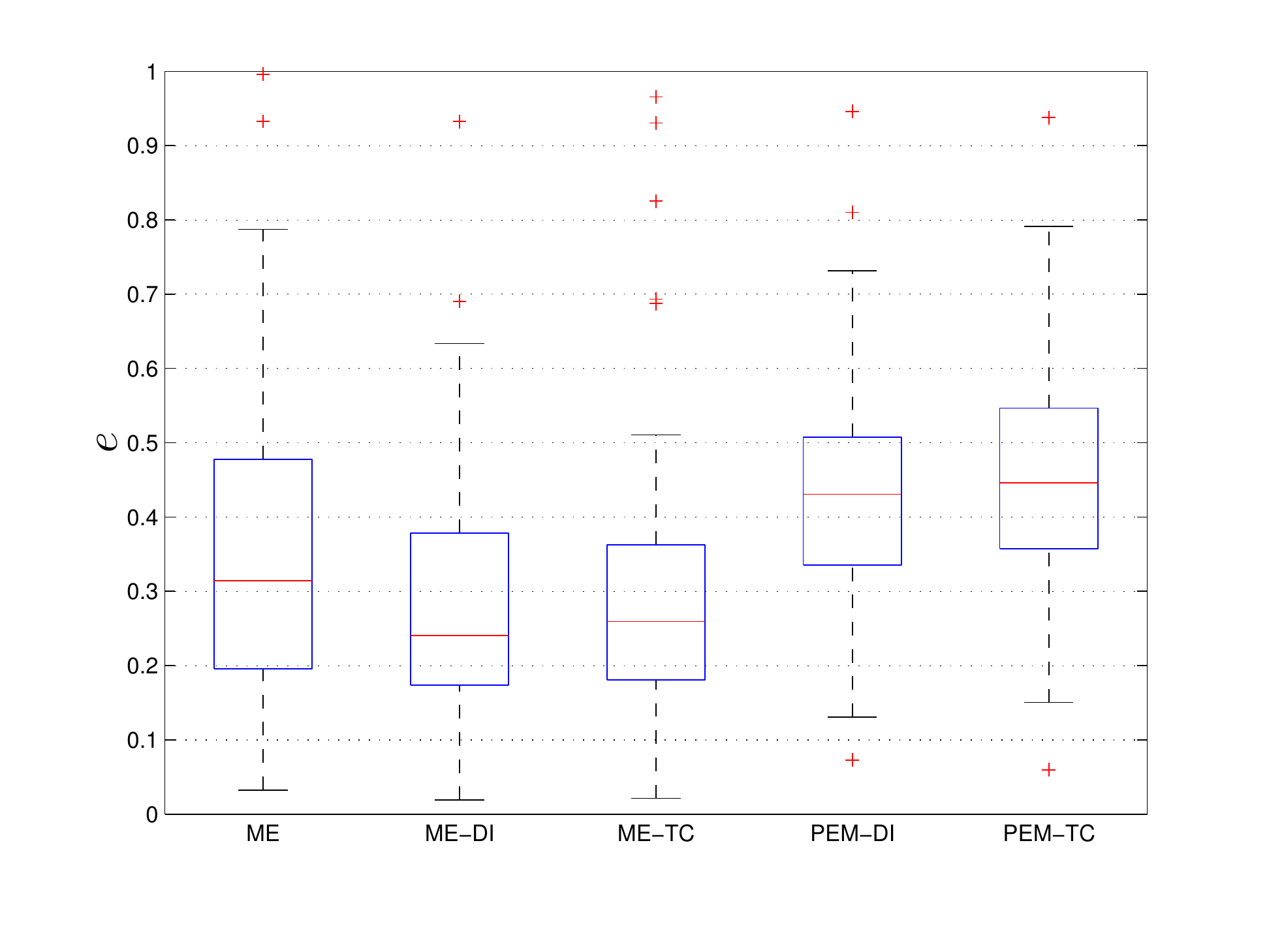}
\caption{Boxplots of the reconstruction error of the spectral estimators in the Monte Carlo experiment. Red crosses represent the outliers.}\label{fig_trial2}
\end{figure} Figure \ref{fig_trial2} shows the boxplots of the reconstruction error. As we can see, ME-DI and ME-TC provide the best performance while the PEM methods provides the worst performance. Finally, in all the runs the kernel based PEM methods returned $\hat b(z)$ with zeros inside the unit circle, so that we did not need to use the refined procedure in \cite{pillonetto2018identification} to approximate $\hat b(z)$ with a minimum phase filter. However choosing the poles of $w(z)$ with absolute value equal to $0.995$, the kernel based PEM methods sometimes returned $\hat b(z)$ non-minimum phase while our methods always returned $\hat b(z)$ minimum phase.

\section{Conclusion}\label{sec:concl}
 In this paper we have introduced a kernel based ME spectral estimator. We have proposed two kernel matrices embedding the property that the impulse response of the one-step ahead predictor should have: regularity and exponential decay to zero. 
Moreover, we showed that the estimated spectrum is given by solving a generalized version of the Yule-Walker equation which leads also a minimum phase spectral factor of the process.   
Finally, numerical experiments showed the effectiveness of the proposed method in respect to the ME spectral estimator and kernel based PEM methods.

\end{document}